\documentclass[11pt]{amsart}
\numberwithin{equation}{section}
\newtheorem{theorem}{Theorem}

\newtheorem{proposition}{Proposition}
\newtheorem{lemma}{Lemma}

\numberwithin{theorem}{section} \numberwithin{lemma}{section}
 \numberwithin{definition}{section}
\numberwithin{proposition}{section}

\def\al{\aligned}
\def\eal{\endaligned}
\def\M{{\bf M}}
\def\be{\begin{equation}}
\def\ee{\end{equation}}
\def\lab{\label}
\def\a{\alpha}
\def\b{\beta}
\def\e{\epsilon}

\def\M{{\bf M}}

\def\al{\aligned}

\def\g{\bar}
\def\p{\partial}
\def\d{\nabla}
\numberwithin{equation}{section}

\begin{document}

\tracingpages 1
\title[]{\bf Isoperimetric inequality under K\"ahler Ricci flow}
\author{Gang Tian and Qi S. Zhang}
\address{BICMR, Peking University, Beijing, 100871, China and Department of Mathematics, Princeton University, Princeton, NJ 02139, USA}
\address{Department of
Mathematics Nanjing University, Nanjing 210093, China; Department
of Mathematics, University of California, Riverside, CA 92521,
USA}
\date{}

\begin{abstract}
Let $({\M}, g(t))$ be a K\"ahler Ricci flow with positive first Chern class. First, we prove a uniform isoperimetric inequality for all time.
 In the process, we also prove a Cheng-Yau type
log gradient bound for positive harmonic functions on $({\M}, g(t))$ and a Poincar\'e inequality without assuming the Ricci
curvature is bounded from below.
\end{abstract}
\maketitle
\section{Introduction}

In this paper we study   K\"ahler Ricci flows
\be
\lab{krf}
\partial_t g_{i\g j} = -  R_{i\g j} + g_{i\g{j}} = \p_i \p_{\g j} u, \quad t>0,
\ee on a compact, K\"ahler manifold $\M$ of complex dimension $m=n/2$, with positive first Chern class.

Given initial K\"ahler metric $g_{i\g j}(0)$, H. D. Cao \cite{Ca:1} proved that (\ref{krf}) has a solution for all time $t$. Recently, many  results
concerning long time and uniform behavior of (\ref{krf}) have appeared. For example, when the curvature
operator or the bisectional curvature is nonnegative, it is known that solutions to (\ref{krf}) stays smooth when
time goes to infinity (see \cite{CCZ:1}, \cite{CT:1} and \cite{CT:2} for examples).
In the general case, Perelman
(cf \cite{ST:1}) proved that the scalar curvature $R$ is uniformly bounded, and the Ricci potential
$u(\cdot, t)$ is uniformly bounded in $C^1$ norm, with respect to $g(t)$.
When the complex dimension $m=2$, let $({\M}, g(t))$ be a solution to (\ref{krf}),  it is proved in
(\cite{CW:1}) that the isoperimetric constant for $({\M}, g(t))$ is bounded
from below by a uniform constant. We mention that an isoperimetric estimate for the Ricci flow on the two sphere was already proven by Hamilton in
\cite{Ha:1}.

In this paper, we prove that in all complex dimensions, the isoperimetric constant for $({\M}, g(t))$ is bounded
from below by a uniform constant. This extends the result of Chen-Wang mentioned above.
 This result seem to add more
weight to
the belief that the K\"ahler Ricci flow converges to a K\"ahler Ricci soliton as $t \to \infty$, except on a
subvariety of complex codimension $2$,

To make the statement precise, let's introduce notations and definition.
We use $\M$ to denote a compact Riemann manifold and $g(t)$ to denote the
metric at time $t$; $d(x, y, t)$ is the geodesic distance under $g(t)$;
$B(x, r, t) = \{ y \in {\M} \ | \ d(x, y, t) < r \}$  is the geodesic ball of radius $r$, under metric $g(t)$,
centered at $x$, and $|B(x, r, t)|_{g(t)}$ is the volume of
$B(x, r, t)$ under $g(t)$;  $d g(t)$ is the volume element. We also
reserve $R=R(x, t)$ as the scalar curvature under $g(t)$. When the time variable $t$ is
not explicitly used, we may also suppress it in the notations mentioned above.

\medskip

The main result of the paper is the following theorem.

\begin{theorem}
\lab{thm1.1}

 Let $({\M}, g(t))$,  $\partial_t g_{i\g j} = -  R_{i\g j} + g_{i\g{j}}$, be a K\"ahler Ricci flow on a $n$ real
 dimensional compact, K\"ahler manifold with positive first Chern class. Then there exists a uniform constant $S_0$,
 depending only on the initial metric $g(0)$ and a numerical constant $C$, such that
 \[
 \left[ \int_{\M} | u |^{n/(n-1)} d g(t) \right]^{(n-1)/n} \le S_0 \int_{\M} | \nabla u |
 dg(t) + \frac{C}{|\M|^{1/n}_{g(t)}} \int_{\M} | u| dg(t)
\]for all $ u \in C^\infty(\M)$.
\end{theorem}

{\it Remark.

 It is well known that Theorem \ref{thm1.1} implies a uniform lower bound for the isoperimetric
constant of $({\M}, g(t))$, i.e. there exists a positive constant $c_0$, depending only on the initial metric such that
\[
I({\M}, g(t)) \equiv \inf_{D \subset \M} \, \frac{ | \partial D |}{ \left[ \min \{ |D|, \, |{\M}-D|\} \right]^{(n-1)/n} } \ge c_0.
\] Here all the volume are with respect to $g(t)$; and $D$ is a subdomain of $\M$ such that
$\p D$ is a $n-1$ dimensional submanifold of $\M$.  A proof can be found in  \cite{CLN:1}  Section 5.1 e.g.

}

\medskip

The proof of the theorem is based on the following properties for  K\"ahler Ricci flow on a compact manifold with positive first Chern class.

\vskip 0.5cm

\noindent {\it Property A.  Let $({\M}, g(t))$ be  a K\"ahler Ricci flow (\ref{krf}) on a compact manifold with positive first Chern class. There exist uniform positive constants $C$, and $\kappa$ so that \\

1. \quad $| R(g(t))| \le C,$

2. $diam (\M, g(t)) \le C,$

3. $\Vert u \Vert_{C^1} \le C.$

4. $|B(x, r, t)|_{g(t)} \ge \kappa r^n$, for all $t>0$ and $r \in (0, diam (\M, g(t)))$. }

5. $
|B(x, r, t)|_{g(t)} \le \kappa^{-1}  r^n
$ for all $r>0$, $t>0$.

\vskip 0.5cm

\noindent {\it Property B. Under the same assumption as in Property A, there exists a uniform constant $S_2$ so that the following $L^2$ Sobolev inequality holds: \\
 \[
 \left( \int_{\M}  v^{2n/(n-2)} d g(t) \right)^{(n-2)/n} \le S_2 \left( \int_{\M} | \nabla v |^2
 dg(t) +  \int_{\M} v^2
 dg(t) \right)
\]for all $ v \in C^\infty(\M, g(t))$.\\
}

Property A  1-4 is due to Perelman (c.f. \cite{ST:1}), Property B was first established in \cite{Z07:1} (see also \cite{Ye:1}, \cite{Z10:1} ).
Property A 5 can be found in \cite{Z11:1} and also \cite{CW:2}.

The rest of the paper is organized as follows. In Section 2, we prove some gradient bounds for
 harmonic functions on $({\M}, g(t))$. Since the bounds do not rely on the usual lower bound of
Ricci curvature, the result may be of independent interest.  Using these bounds, we prove the
theorem in Section 3.

\medskip

\section{gradient bounds for harmonic functions}

In order to prove the theorem, in this section we state and prove a number of results on harmonic functions on certain manifolds with fixed metric.
These results are well known when
the manifold has nonnegative Ricci curvature, a property that is not available for us.
Since some of these results may be of independent interest, we will also deal with the real
variable case and impose some conditions which are more general than needed for the
proof of the theorems in Section 1.  As the metric is independent of time in this section, we will
suppress the time variable $t$.

In this section, the basic assumptions on the $n$ real dimensional manifolds $\M$ are

{\it Assumption  1. $L^2$ Sobolev inequality: there is a positive constant $\a$ such that
\[
 \left( \int_{\M}  u^{2n/(n-2)} d g(t) \right)^{(n-2)/n} \le \a \left( \int_{\M} | \nabla u |^2
 dg(t) +  \int_{\M} u^2
 dg(t) \right)
\]for all $ u \in C^\infty(\M)$.
}

{\it Assumption  2. There exists a positive constant $\kappa$, such that
\[
\kappa r^n \le |B(x, r)| \le \kappa^{-1} r^n, \qquad x \in \M, \quad 0 < r< diam(\M) \le 1.
\]}

{\it Assumption 3. There exists a smooth function $L=L(x)$ and 2 smooth parallel $(2, 2)$ tensor fields $P$
 and $Q$ such that the Ricci curvature is given by
 \[
 R_{ij} =P^{kl}_{ij} \p_k \p_l L + Q^{kl}_{ij} g_{kl}
 \] under a local coordinates. Moreover $\Vert P \Vert_\infty \le 1$, $\Vert Q \Vert_\infty \le 1$. Here $\p_k \p_l L $ is the
 Hessian of $L$.
 }

Note that Assumption 3 includes  as a special case, the formula for the Ricci curvature on K\"ahler manifolds
$\p_i\p_{\g j} u = g_{i\g j}-R_{i\g j}$ where $u$ is the Ricci potential.

\begin{lemma}
\lab{ledumean}
Suppose $(\M, g)$ is a compact Riemann manifold of real dimension $n$,  satisfying
Assumptions 1, 2, 3.

 Let $u$ be a smooth harmonic function in $B(x_0, r)$ where $x_0 \in \M$ and
 $ r \le diam(\M)$.  Then there exists a positive constant $C_0=C_0(\a, \kappa, \Vert \d L \Vert_\infty)$
 such that
 \[
 \sup_{x \in B(x_0, r/2)} |\nabla u(x)|  \le C_0 \frac{1}{r}
 \left(\frac{1}{r^n} \int_{B(x_0, r)} u^2 dg \right)^{1/2}.
 \]
\end{lemma}

\proof  Since $u$ solves $\Delta u=0$, by Bochner's formula, we have
\be
\lab{dddu}
\Delta |\nabla u |^2 = 2 | Hess \, u |^2 + 2 Ric (\nabla u, \nabla u).
\ee  Given $\sigma \in (0, 1)$, let $\psi=\psi(x)$ be a standard Lipschitz cut-off function such that
$\psi(x)=0$ when $x \in B(x_0, r)^c$; $0 \le \psi \le 1$ and $\psi(x)=1$ when $x \in
B(x_0, \sigma r)$ and $|\nabla \psi| \le \frac{4}{(1-\sigma) r}$. We mention that
no second order derivatives of $\psi$ are involved. So we only need $\psi$ is
Lipschitz.

For clarity of presentation, we write
\[
f = |\nabla u|^2.
\]Given a number $p \ge 1$, using $f \psi^2$ as a test function on (\ref{dddu}), after a routine calculation, we derive
\be
\lab{dfpsi}
\al
&\frac{2p-1}{p^2} \int_{B(x_0, r)} |\nabla (f^p \psi)|^2 dg \\
&\le \frac{C}{(1-\sigma)^2 r^2} \int_{B(x_0, r)} f^2 dg
-2 \int_{B(x_0, r)}  | Hess \, u |^2 f^{2p-1} dg -
2 \int_{B(x_0, r)} Ric (\nabla u, \nabla u) f^{2p-1} \psi^2 dg\\
&\equiv I_1 + I_2 +I_3.
\eal
\ee Now we want to absorb part of $I_3$ by $I_2$ which is a good term.
In a local orthonormal coordinates, we denote $u_i$ the i-th component of $\nabla u$. Then
By Assumption 3, we have, after integrating by parts,
\[
\al
I_3 &= - 2 \int R_{ij} u_i u_j f^{2p-1} \psi^2 dg\\
&=-2 \int P^{kl}_{ij} (\p_k \p_l L)  \, u_i u_j f^{2p-1} \psi^2 dg - 2 \int Q^{kl}_{ij} g_{kl} u_i u_j f^{2p-1} \psi^2 dg\\
&=2 \int P^{kl}_{ij} (\p_l L)  \, (\p_k u_i) u_j f^{2p-1} \psi^2 dg
+ 2 \int P^{kl}_{ij} (\p_l L)  \, u_i (\p_k u_j) f^{2p-1} \psi^2 dg\\
& \qquad+2 \int P^{kl}_{ij} (\p_l L)  \, u_i  u_j \p_k (f^{2p-1} \psi^2) dg
- 2 \int Q^{kl}_{ij} g_{kl} u_i u_j f^{2p-1} \psi^2 dg.
\eal
\] Here we also used the assumption that the $P$ tensor is parallel.
To control the second from last term in the above identity, we notice that
$
|u_i u_j | \le |\d u|^2 = f$ and that
\[
f \p_k (f^{2p-1} \psi^2) =(\p_k f^p) \psi^2 \frac{2p-1}{p} f^p
+ 2 f^{2p} (\p_k \psi) \psi.
\]
From here, using Young's
inequality, it is easy to see that
\[
I_3 \le \frac{1}{2} \int |\nabla (f \psi)|^2 dg + \int | Hess \, u |^2 f^{2p-1} \psi^2 dg
+ C \left(\frac{ \Vert \d L \Vert^2_\infty}{[(1-\sigma) r]^2} + \Vert \d L \Vert^2_\infty +1 \right)
\int f^{2p} \psi^2 dg.
\]Substituting this to (\ref{dfpsi}), we deduce
\[
\int |\nabla (f^p \psi)|^2 dg \le C p \left(\frac{ \Vert \d L \Vert^2_\infty}{[(1-\sigma) r]^2} + \Vert \d L \Vert^2_\infty +1 \right)
\int f^{2p} \psi^2 dg.
\] Since $diam (\M) \le 1$ by Assumption 2, the last inequality implies
\[
\int |\nabla (f^p \psi)|^2 dg \le C p \frac{ \Vert \d L \Vert^2_\infty + 1}{[(1-\sigma) r]^2}
\int f^{2p} \psi^2 dg.
\] Using the $L^2$ Sobolev inequality in Assumption 1 and Moser's iteration, we deduce, by choosing $p= (n/(n-2))^i$, and replacing $(1-\sigma)r$ by
$\left(1/4\right)^i r$,  $i=0, 1, 2, ...$,
that
\be
\lab{du3r/4}
 \sup_{x \in B(x_0, r/2)} |\nabla u(x)|^2  \le C_1
 \frac{1}{r^n} \int_{B(x_0, 3r/4)} |\nabla u|^2 dg
 \ee where $C_1=C_1(\a, \kappa, \Vert \d L \Vert_\infty)$. We observe that even
though that the number $p$ appears on the right hand side of the inequality
before (\ref{du3r/4}), but its growth is only an exponential of $i$.  As well known,  it will be suppressed by the Moser iteration process, just like the term $1/(1-\sigma)^2$.

 Next we take $\sigma=3/4$ in the definition of the cut off function $\psi$. Using $u \psi^2$
 as a test function on $\Delta u =0$,  we infer, after a routine calculation
 \[
  \int_{B(x_0, 3r/4)} |\nabla u|^2 dg \le \frac{C}{r^2} \int_{B(x_0, r)} u^2 dg.
 \] Here $C$ is a numerical constant. Combining the last two inequalities we arrive at
 \[
 \sup_{x \in B(x_0, r/2)} |\nabla u(x)|^2  \le C_0
 \frac{1}{r^{n+2}} \int_{B(x_0, r)} u^2 dg
 \]
  where $C_0=C_0(\a, \kappa, \Vert \nabla L \Vert_\infty)$.
 \qed

 The next lemma is simply the $L^2$ mean value inequality for the Laplace and heat equation
 under Assumptions 1 and 2. Since the result is well known (Grigoryan \cite{Gr:1} and Saloff-Coste \cite{Sa:1}), we
 omit the proof.

 \begin{lemma}
 \lab{lemvp}
  Let $\M$ be a manifold satisfying  Assumptions 1, 2.

  Suppose $u$ be is a smooth harmonic function in $B(x_0, r)$ where $x_0 \in \M$ and
 $ r \le diam(\M)$.  Then there exists a positive constant $C_1=C_1(\a, \kappa)$
 such that
 \[
 \sup_{x \in B(x_0, r/2)} |u(x)|  \le C_1
 \left(\frac{1}{r^n} \int_{B(x_0, r)} u^2 dg \right)^{1/2}.
 \]

 Suppose $u$ is a solution of the heat equation $\Delta u -\partial_t u =0$ in the space time
 cube $B(x_0, r) \times [t_0-r^2, t_0]$. Then
  \[
 \sup_{(x, t) \in B(x_0, r/2) \times [t_0 - r^2/4]} |u(x, t)|  \le C_1
 \left(\frac{1}{r^{n+2}} \int^t_{t-r^2}\int_{B(x_0, r)} u^2 dg ds \right)^{1/2}.
 \]
\end{lemma}

The next lemma provides bounds for the Green's function of the Laplacian and its gradients.

 \begin{lemma}
 \lab{leDGbound}
 Let $\M$ be a manifold satisfying  Assumptions 1, 2 and 3.
 Assume also $diam(\M)>\beta>0$ for a positive constant $\beta$.
 Let $\Gamma_0$ be the Green's function of the Laplacian $\Delta$ on $\M$.
 Then there exists
 a positive constant $C_0=(\a, \beta, \kappa, \Vert \nabla L \Vert_\infty)$ such that

 (a). $ |\Gamma_0(x, y)| \le \frac{C_0}{d(x, y)^{n-2}}$, $x, y \in \M$,

 (b). $ |\nabla_x \Gamma_0(x, y)| \le \frac{C_0}{d(x, y)^{n-1}}$, $x, y \in \M$.
\end{lemma}

\proof Once (a) is proven, (b) is just a consequence of (a) and Lemma \ref{ledumean} applied
on the ball $B(x, d(x, y)/2)$.  So now we just need to prove (a).

On a compact manifold $\M$, we know that \be \lab{gammaG}
\Gamma_0(x, y) = \int^\infty_0 \left( G(x, t, y)- \frac{1}{|\M|}
\right) dt \ee where $G$ is the fundamental solution of the heat
equation $\Delta u - \partial_t u=0$. We remark that the metric is
fixed here. So we need to  bound  $G$. Under Assumptions 1 and 2,
Grigoryan \cite{Gr:1} and Saloff-Coste \cite{Sa:1} proved that
there exist positive constants $A_1, A_2, A_3$ which depend only
on $\a$ and $\kappa$ such that
\be
\lab{Gaub} G(x, t, y) \le A_1 (
1 + \frac{1}{t^{n/2}} ) e^{- A_2 d(x, y)^2/t}.
\ee Fixing $x, y$ and $t$, we write $u=u(z, l) =G(z, l, y)$ and regard it as a solution of the
heat equation in the cube $B(x, r) \times [t-r^2, t]$. Here $r=\sqrt{t}/2$.
Extending Lemma \ref{ledumean} to the parabolic case in a routine manner, we know that
\[
 |\nabla u(x, t)|  \le \frac{C_1}{r}
 \left(\frac{1}{r^{n+2}} \int^t_{t-r^2}\int_{B(x, r)} u^2 dg ds \right)^{1/2}.
\]Substituting (\ref{Gaub}) to the right hand side, we know that
\[
| \nabla_x G(x, t, y)| \le A_1 (
1 + \frac{1}{t^{(n+1)/2}} ) e^{- A_2 d(x, y)^2/t}.
\]Here the constants $A_1$ and $A_2$ may have changed.
It is well known that this gradient bound and the upper bound (\ref{Gaub})
together imply a Gaussian lower bound for the heat kernel $G$. See \cite{CD:1}, p1165 e.g.
Now, by
\cite{Sa:1}, the following $L^2$ Poincar\'e inequality holds: for any $u
\in C^\infty(\M)$, $r \in (0, diam(\M)]$,
\be
\lab{l2Poin}
\int_{B(x, r/2)} | u - \g u_{B(x, r/2)}|^2 dg \le A_3 r^2 \int_{B(x, r)} |\nabla u|^2 dg.
\ee By a trick in Jerison \cite{J:1}, which uses only volume doubling property, one has that
\be
\lab{goodl2Poin}
\int_{B(x, r)} | u - \g u_{B(x, r)}|^2 dg \le C A_3 r^2 \int_{B(x, r)} |\nabla u|^2 dg.
\ee Here $C$ depends only on $\kappa$. We mention that some of the cited results
were stated for complete noncompact manifolds. But they are also valid for complete,
closed manifolds as long as the diameters are uniformly bounded.

Let $u_0 \in C^\infty(\M)$ be a function such that $\int_\M u_0 dg =0$. Then the function
\be
\lab{uxt=}
u(x, t) = \int_\M \left( G(x, t, z)- \frac{1}{|\M|} \right) u_0(z) dg(z)
\ee is a solution to the heat equation such that $\int_\M u(x, t) dg(x)=0$.
By the $L^2$ Poincar\'e inequality with $r=diam (\M)$, we have
\[
\int_\M u^2 dg \le C \, A_3 \, diam(\M)^2 \int_\M |\nabla u|^2 dg \le  C A_3  \int_\M |\nabla u|^2 dg
\]since $diam(\M) \le 1$ by assumption.  From this we deduce
\[
\frac{d}{d t} \int_\M u^2 dg = - 2 \int_\M |\nabla u|^2 dg = -2 (CA_3)^{-1}  \int_\M u^2 dg
\] and consequently
\[
\int_\M u^2(z, s) dg \le e^{ - 2 (CA_3)^{-1} s}  \int_\M u^2_0(z) dg, \qquad s>0.
\]

Recall that we assume $diam(\M)>\beta>0$. For $t \ge \beta^2$, we can apply Lemma
\ref{lemvp} to get
\[
u^2(x, t) \le  C^2_1
\frac{1}{\beta^{n+2}}  \int^t_{t-\beta^2} \int_{\M} u^2(z, s) dg ds.
\] Combining this with the previous inequality, we arrive at
\[
u^2(x, t) \le  C_2 e^{ - 2 (CA_3)^{-1} t}  \int_\M u^2_0(z) dg
\] where $C_2=C_0(\a, \b, \kappa, A_3)$.
By (\ref{uxt=}), this means
\[
\left[ \int_\M \left( G(x, t, z)- \frac{1}{|\M|} \right) u_0(z) dg \right]^2
\le C_2 e^{ - 2 (CA_3)^{-1} t} \int_\M u^2_0(z) dg
\] Fixing $x \in \M$ and $t \ge \beta^2$, and taking $u_0(z) = G(x, t, z)- \frac{1}{|\M|}$ in the above
inequality, we obtain
\be
\lab{intg-1m}
 \int_\M \left( G(x, t, z)- \frac{1}{|\M|} \right)^2 dg \le C_2 e^{ - 2 (CA_3)^{-1} t},
 \quad t \ge \b^2.
 \ee Fixing $x$, the function $h(z, t) \equiv G(x, t, z)- \frac{1}{|\M|}$ is also a solution to the heat
 equation. Applying the mean value inequality in Lemma \ref{lemvp} on the cube
 $B(y, \beta) \times [t-\beta^2, t]$, we infer
\[
h^2(y, t) \le  C^2_1
\frac{1}{\beta^{n+2}}  \int^t_{t-\beta^2} \int_{\M} h^2(z, s) dg ds.
\] That is
\[
\left(G(x, t, y)- \frac{1}{|\M|}\right)^2 \le C^2_1
\frac{1}{\beta^{n+2}}  \int^t_{t-\beta^2} \int_{\M}  \left( G(x, s, z)- \frac{1}{|\M|} \right)^2 dg ds.
\]Substituting (\ref{intg-1m}) to the last inequality, we deduce
\be
\lab{Gboundt>b}
| G(x, t, y)- \frac{1}{|\M|} | \le C_3 e^{- C_4 t}, \qquad t \ge \b^2,
\ee where $C_3, C_4$ depend only on $\a, \b, \kappa$ and $A_3$ which only depends on $\a,
\kappa$.

From (\ref{gammaG}),
\[
\al
\Gamma_0(x, y) &= \int^{\b^2}_0  \left( G(x, t, y)- \frac{1}{|\M|} \right) dt +
\int^\infty_{\b^2}  \left( G(x, t, y)- \frac{1}{|\M|} \right) dt\\
&\equiv I_1 + I_2.
\eal
\]Using the bound (\ref{Gaub}) on $I_1$ and (\ref{Gboundt>b}) on $I_2$, we derive,
after simple integration,
\[
|\Gamma_0(x, y) | \le \frac{C_0}{d(x, y)^{n-2}},
\]where $C_0$ depends only on  $\a, \b, \kappa$ .  This proves part (a) of the Lemma.

As mentioned earlier, part (b) follows from part (a) and Lemma \ref{ledumean}.
\qed

The next result is a Cheng-Yau type log gradient estimate. Although not used in the proof
of the theorems, it may be of independent interest.

\begin{proposition}
Let $\M$ be a manifold satisfying  Assumptions 1, 2 and 3.  Let $u$ be a positive harmonic
function in the geodesic ball $B(x, 2r)$, which is properly contained in $\M$. Then there
exists a positive constant $C$, depending only on the controlling constants in Assumptions
1-3, such that
\[
\sup_{B(x, r)} | \nabla \ln u | \le \frac{C}{r}
\] when $ r \in (0, 1]$.
\proof
\end{proposition}

For convenience, we use the following notations
\[
h \equiv \ln u, \quad F \equiv | \nabla h|^2.
\] Following \cite{CY:1}, it is well known that $\Delta h = - F$ and
\[
\Delta F = - 2 \nabla h \nabla F + 2 | Hess \, h|^2 + 2 Ric (\d h, \d h).
\] Consider the function
\be
\lab{w=}
w  \equiv F^{5n}.
\ee By a routine calculation, we know that, for any $p \ge 1$,
\be
\lab{ddwp}
\Delta w^p  \ge - 2 \nabla h \nabla w^p + 10 n p F^{5np-1} | Hess \, h|^2 + 10 n p F^{5np-1} Ric (\d h, \d h)
\ee

Given $\sigma \in (0, 1)$, let $\psi=\psi(x)$ be a standard smooth cut-off function such that
$\psi(x)=0$ when $x \in B(x_0, r)^c$; $0 \le \psi \le 1$ and $\psi(x)=1$ when $x \in
B(x_0, \sigma r)$ and $|\nabla \psi| \le \frac{4}{(1-\sigma) r}$.
Using $w^2 \psi$ as a test function on (\ref{ddwp}),  we deduce, after a straight forward
calculation, that
\be
\lab{dwpp}
\al
\int |\d(w^p \psi)|^2 dg &\le -10 n p \int F^{5np-1} \, | Hess \, h|^2 w^p \psi^2 dg +
2 \int  \nabla h \nabla w^p \, w^p \psi^2 dg\\
&\qquad  -
10 n p \int F^{5np-1} Ric (\d h, \d h) w^p \psi^2 dg\\
&\equiv I_1 + I_2 + I_3.
\eal
\ee

Next we will show that the negative term $I_1$ dominate $I_2$ and $I_3$, modulo
some harmless terms.
Observe that
\[
\al
I_2 &= \int \psi^2 \d h \d w^{2p} dg \\
&= -2 \int \psi \d \psi \d h \, w^{2p} dg  -
\int \psi^2 \Delta h \, w^{2p} dg.
\eal
\] Recall that $\Delta h = -|\d h|^2 =F $. Hence, by Young's inequality, for any given $\e>0$,
\be
\lab{i2<}
I_2 \le (\e +1) \int F w^{2 p} \psi^2 dg + \e^{-1}  \Vert \d \psi \Vert^2_\infty
\int w^{2p} \psi^2 dg.
\ee

It takes a little longer to prove the bound for $I_3$.
By our condition on the Ricci curvature $R_{ij}$, we have
\[
I_3 = -10 n p \int F^{5np-1} ( P^{kl}_{ij} \p_k \p_l L + Q^{kl}_{ij} g_{kl}) \,
\p_i h \p_j h \, w^p \psi^2 dg.
\] After integration by parts, this becomes
\be
\lab{i3=}
\al
I_3 &= 10 np (5np-1) \int F^{5 np-2} \p_k   F  \, P^{kl}_{ij} \p_l L \, \p_i h \p_j h \, w^p \psi^2 dg \\
& \qquad + 10 n p  \int F^{5 np-1}   \, P^{kl}_{ij} \p_l L \, (\p_k \p_i h) \,  \p_j h \, w^p \psi^2 dg\\
& \qquad + 10 n p  \int F^{5 np-1}   \, P^{kl}_{ij} \p_l L \, \p_i h \, (\p_k \p_j h) \, w^p \psi^2 dg \\
&  \qquad + 10 n p  \int F^{5 np-1}   \, P^{kl}_{ij} \p_l L \, \p_i h \,  \p_j h \,  \p_k (w^p \psi) \psi  dg \\
&  \qquad + 10 n p  \int F^{5 np-1}   \, P^{kl}_{ij} \p_l L \, \p_i h \,  \p_j h \,  w^p \psi  \p_k \psi\\
&\qquad -10 n p \int F^{5np-1}  Q^{kl}_{ij} g_{kl}\,
\p_i h \p_j h \, w^p \psi^2 dg\\
&\equiv T_1 + ... + T_6.
\eal
\ee

Let us bound $T_i$, $i=1, ..., 6$. Observe that
\[
| T_1 | \le 10 np (5np-1) \Vert  \nabla L \Vert_\infty \int F^{5 np-2} |\nabla F |  \, | \nabla h |^2 \, w^p \psi^2 dg.
\] Since $| \nabla h |^2 = F$, we deduce, using $w^p=F^{ 5np}$,
\[
\al
| T_1 | &\le 10 np (5np-1) \Vert \nabla L \Vert_\infty \int F^{5 np-1} |\nabla F |  \, w^p \psi^2 dg\\
&\le 10 np   \Vert \nabla L \Vert_\infty \int  |\nabla w^p |  \, w^p \psi^2 dg.
\eal
\]  Thus, after a little calculation, we obtain,
\be
\lab{t1<}
| T_1 | \le \frac{1}{10} \int | \nabla (w^p \psi) |^2 dg + c p^2  \Vert  \nabla L \Vert^2_\infty \int   w^{2 p} \psi^2 dg + c \Vert \d \psi \Vert^2_\infty \int_{supp \, \psi}  w^{2 p}  dg.
\ee Next
\[
\al
| T_2 | &\le 10 np  \Vert \nabla L \Vert_\infty \int F^{5 np-1} \, | Hess \, h | \, | \nabla h |
\, w^p \psi^2 dg \\
& \le  np \int F^{5 np-1} \, | Hess \, h |^2 \,
\, w^p \psi^2 dg + c np \Vert \nabla L \Vert^2_\infty \int F^{5 np-1} \, |\nabla h |^2 \,
\, w^p \psi^2 dg.
\eal
\] Recalling again that $|\nabla h |^2 = F$ and the definition of $I_1$, we deduce
\be
\lab{t2<}
| T_2 | \le - \frac{I_1}{10} +
c np \Vert  \nabla L \Vert^2_\infty \int
\, w^{2 p} \psi^2 dg.
\ee Since $T_3$ is similar to $T_2$, we also have
\be
\lab{t3<}
| T_3 | \le - \frac{I_1}{10} +
c np \Vert \nabla L \Vert^2_\infty \int
\, w^{2 p} \psi^2 dg.
\ee

By Young's inequality
\[
| T_4 | \le
\frac{1}{2} \int | \nabla ( w^p \psi)|^2 dg + 50 n^2 p^2
 \Vert \nabla L \Vert^2_\infty \int F^{10 np-2} \, |\nabla h |^4  \psi^2\,
 dg.
\] Since $F = |\d h|^2 $ and $w=F^{5n}$, this shows
\be
\lab{t4<}
| T_4 | \le
\frac{1}{2} \int | \nabla ( w^p \psi)|^2 dg + c p^2
 \Vert \nabla L \Vert^2_\infty \int w^{2 p} \psi^2\,
 dg.
\ee Next
\[
|T_5| \le 10 np \Vert \nabla L \Vert_\infty \, \Vert \psi \Vert_\infty
\int F^{5np-1} | \d h|^2 w^p \psi dg,
\] which becomes
\be
\lab{t5<}
|T_5| \le 10 np \Vert \nabla L \Vert_\infty \, \Vert \psi \Vert_\infty
\int w^{2p} \psi dg.
\ee Lastly
\be
\lab{t6<}
| T_6 | \le 10 np
\int F^{5np-1} | \d h|^2 w^p \psi^2 dg = 10 np
\int  w^{2 p} \psi^2 dg.
\ee Substituting (\ref{t1<})-(\ref{t6<}) into (\ref{i3=}), we find that
\be
\lab{i3<}
| I_3 | \le \frac{|I_1|}{5} + \frac{3}{5} \int | \nabla ( w^p \psi)|^2 dg
+ c \frac{ p^2 \Vert \nabla L \Vert^2_\infty +1}{[(1-\sigma) r]^2} \int_{supp \, \psi}
w^{2p} dg.
\ee
 Here we recall that
\[
I_1=
 -10 n p \int F^{5np-1} \, | Hess \, h|^2 w^p \psi^2 dg.
 \]Using the inequality
 \[
 | Hess \, h |^2 \ge \frac{1}{n} ( \Delta h )^2 = \frac{1}{n} |\d h |^4,
 \]we find that
 \[
 I_1 = \frac{I_1}{2} +  \frac{I_1}{2} \le  \frac{I_1}{2} - 5p \int
F^{5np-1} \, | \d h|^4 w^p \psi^2 dg
\] which induces, since $w=F^{ 5n}$ and $F = | \d h |^2$, that
\be
\lab{i1<2}
I_1  \le \frac{I_1}{2} - 5p \int
F w^{2 p} \psi^2 dg.
\ee Substituting  (\ref{i1<2}),  (\ref{i3<}) and  (\ref{i2<}) into (\ref{dwpp}), we deduce
\[
\al
\int |\d(w^p \psi)|^2 dg &\le
\frac{I_1}{2} - 5p \int
F w^{2 p} \psi^2 dg +
 (\e +1) \int F w^{2 p} \psi^2 dg + \e^{-1}  \Vert \d \psi \Vert^2_\infty
\int w^{2p} \psi^2 dg \\
&\qquad + \frac{|I_1|}{5} + \frac{3}{5} \int | \nabla ( w^p \psi)|^2 dg
+ c \frac{ p^2 \Vert \nabla L \Vert^2_\infty +1}{[(1-\sigma) r]^2} \int_{supp \, \psi}
w^{2p} dg.
\eal
\] Since $p \ge 1$, we can take$\e=1$ and obtain
\[
\int |\d(w^p \psi)|^2 dg + \int (w^p \psi)^2 dg \le c \frac{ p^2 \Vert \nabla L \Vert^2_\infty +1}{[(1-\sigma) r]^2} \int_{supp \, \psi}
w^{2p} dg
\]where $c$ may have changed in value. By the Sobolev inequality in Assumption 1,
this implies
\[
 \left( \int (w^p \psi)^{2n/(n-2)} dg \right)^{(n-2)/n}  \le c \a \frac{ p^2 \Vert
\nabla L \Vert^2_\infty +1}{[(1-\sigma) r]^2} \int_{supp \, \psi}
w^{2p} dg.
\] From this, the standard Moser's iteration implies
\[
\sup_{B(x, \sigma r)} w^2 \le \frac{C(\a, n,  \Vert  \nabla L \Vert^2_\infty)}{(1-\sigma)^n r^n} \int_{B(x, r)} w^2 dg
\]for $r, \sigma \in (0, 1]$.
Using $w=F^{5n}$, we arrive at
\[
\sup_{B(x, \sigma r)} F \le \left( \frac{C(\a, n,  \Vert  \nabla L \Vert^2_\infty)}{(1-\sigma)^n r^n} \int_{B(x, r)} F^{10n} dg \right)^{1/(10n)}
\]for $r, \sigma \in (0, 1]$. Using the volume doubling property and an  algebraic trick in
\cite{LS:1} e.g., we deduce
\be
\lab{F<}
\sup_{B(x, r/2)} F \le \frac{C(\a, n,  \Vert \nabla L \Vert^2_\infty)}{r^n} \int_{B(x, r)} F dg
\ee for $r, \sigma \in (0, 1]$. Using integration by parts, it is known that
\[
 \int_{B(x, r)} F dg =  \int_{B(x, r)} |\nabla (\ln u)|^2 dg \le 4 \frac{|B(x, 4 r)|}{ r^2} \le c r^{n-2}
\]where we have used Assumption 2. Substituting this to (\ref{F<}), we arrive at
\[
\sup_{B(x, r/2)} |\nabla (\ln u)| \le \frac{C(\a, n,  \Vert \nabla L \Vert^2_\infty)}{r}
\]proving the proposition. \qed

\section{Proof of the Theorem}

\proof (Theorem \ref{thm1.1}).

For simplicity of presentation, we omit the time variable in the proof. It is also clear that we
can take $\g u =0$.

{\it Step 1.}
\medskip

Pick $ u \in C^\infty(\M)$. Since $\Delta u = \Delta u$ and $\g u=0$, we have
\[
u(x) = -\int_\M \Gamma_0(x, y) \Delta u(y) dg(y),
\] where $\Gamma_0$ is the Green's function of the Laplacian on $\M$.
Pick a small balls $B(x, r)$. Then,
\[
\al
u(x) &= -\lim_{r \to 0} \int_{\M-B(x, r)} \Gamma_0(x, y) \Delta u(y) dg(y)\\
& = \lim_{r \to 0} \int_{\M-B(x, r)} \nabla \Gamma_0(x, y) \nabla u(y) dg(y)
- \lim_{r \to 0}\int_{\partial B(x, r)}  \Gamma_0(x, y) \partial_n  u(y) dS.
\eal
\] Here we have used
 integration by parts.  Note that $|\Gamma_0(x, y)| \le \frac{C_0}{d(x, y)^{n-2}}$
by Lemma\ref{leDGbound}.  Also the volume of $\partial B(x, r)$, the small spheres of radius $r$,
is bounded from above by $C r^{n-1}$. So the second limit is $0$. We mention
that one does not need a uniform in time bound for $|\partial B(x, r)|$ since
we are freezing a time and taking the limit $r \to 0$.
Hence
\[
u(x) = \int_\M \nabla \Gamma_0(x, y) \nabla u(y) dg(y).
\] According to Lemma \ref{leDGbound}, this implies
\be
\lab{u<I1du}
|u(x)| \le C_0 \int_\M  \frac{|\nabla u(y)|}{d(x, y)^{n-1}} dg(y) \equiv C_0 I_1(|\nabla u|) (x).
\ee Here $I_1$ is the Riesz potential of order $1$.

We claim that there exists a constant $C_1$, depending only on the
constant $\kappa$ in Property A 5, such that \be \lab{I1mf}
|I_1(f)(x)| \le C_1 [ M(f) (x)]^{1-(1/n)} \, \Vert f
\Vert^{1/n}_1. \ee for all smooth function $f$ on $\M$.
Here $M(f)$ is the Hardy-Littlewood maximal function.
The proof
given here is more or less the same as in the Euclidean case (p86
\cite{Zi:1}), under Property A 5, i.e.  $\kappa r^n \le |B(x, r)|
\le \kappa^{-1} r^n$.   Let $\delta$ be a positive number, then
\[
\al
|I_1(f)(x)| &\le \int_{B(x, \delta)} \frac{|f(y)|}{d(x, y)^{n-1}} dg +
\int_{B^c(x, \delta)} \frac{|f(y)|}{d(x, y)^{n-1}} dg\\
&\le \Sigma^\infty_{j=0} \int_{\{ 2^{-j-1} \delta \le d(x, y) <2^{-j} \delta \}}
\frac{|f(y)|}{d(x, y)^{n-1}} dg + \delta^{1-n} \int_\M |f(y)| dg\\
& \le \Sigma^\infty_{j=0}  (2^{(j+1)}/\delta)^{n-1} |B(x, 2^{-j} \delta)|
\frac{1}{|B(x, 2^{-j} \delta)|}  \int_{B(x, 2^{-j} \delta)} |f(y)|dg + \delta^{1-n} \int_\M |f(y)| dg\\
&\le \Sigma^\infty_{j=0}  (2^{(j+1)}/\delta)^{n-1} |B(x, 2^{-j} \delta)|
 \, M(f)(x) + \delta^{1-n} \Vert f \Vert_1.
\eal
\] By Property A 5,
\[
|B(x, 2^{-j} \delta)| \le \kappa^{-1} (2^{-j} \delta)^n.
\] Combining the last 2 inequalities we deduce
\[
|I_1(f)(x)| \le  C \kappa^{-1} \delta \, M(f)(x) + \delta^{1-n} \Vert f \Vert_1
\] which implies (\ref{I1mf}) by taking $\delta =[ M(f)(x)/\Vert f \Vert_1]^{-1/n}$.
We remark that if $\delta>diam (\M)$, then the integral $\int_{B^c(x, \delta)} \frac{|f(y)|}{d(x, y)^{n-1}} dg$ is regarded as zero.

Since Property A 4-5 induces volume doubling property, it is well known that
the maximal operator is bounded from $L^1(\M)$ to weak $L^1(\M)$, i.e.
there is a positive constant $C_2$, depending only on $\kappa$ such that
\[
\b | \{ x \, | \, M(f)(x)> \b \} | \le C_2  \Vert f \Vert_1, \qquad
\]for all $\b>0$. A short proof can be found
in Chapter 3 of Folland's book \cite{Fo:1} e.g. Note the proof there is written
for the Euclidean space. But as indicated below,  it is clear that it works for all metric spaces with
volume doubling property.  Pick $x \in S_\b \equiv \{ x \, | \, M(f)(x)> \b \}$. Then by definition of $M(f)(x)$, there
exists radius $r_x>0$ such that
\[
\frac{1}{|B(x, r_x)|} \int_{B(x, r_x)}  |f(y)| dg > \b.
\]Note that the family of balls $\{ B(x, r_x) \, |  \, x \in S_\b \}$ is an
open cover of $S_\b$. Since the manifold is compact, by well known covering
argument for compact metric spaces,  there exists a finite subfamily $\{ B(x, r_{x_i}) \, |  \, i=1, ..., m \} $
of disjoint balls
such that $\{ B(x, 3 r_{x_i}) \, |  \, i=1, ..., m \}$ covers $S_\b$. Using
  volume doubling property, one has
\[
\b |S_\b| \le \b \Sigma_i  |B(x, 3 r_{x_i})| \le C \Sigma_i  \b |B(x,  r_{x_i})|
\le  \Vert f \Vert_1.
\]
Combining this with (\ref{I1mf}), we obtain, for all $\a>0$,
\[
\al
 | \{ x \, | \, I_1(f)(x)> \a \} |& \le  | \{ x \, | \, M(f)(x)> \frac{\a^{n/(n-1)}}{\Vert f \Vert^{1/(n-1)}_1
 C^{n/(n-1)}_1}\} | \\
 &\le C_2 C^{n/(n-1)}_1 \Vert f \Vert^{1/(n-1)}_1 \a^{-n/(n-1)}  \Vert f \Vert_1.
 \eal
 \]Thus
 \be
 \lab{I1weak}
\a^{n/(n-1)}  | \{ x \, | \, I_1(f)(x)> \a \} | \le C_2 C^{n/(n-1)}_1 \Vert f \Vert^{n/(n-1)}_1
\ee By (\ref{u<I1du}) we have
\[
| \{ x \, | \, |u(x)|>\a \} \le | \{ x \, | \, |I_1(\nabla u)(x)|>\a C^{-1}_0 \} ,
\] which infers, via (\ref{I1weak}) with $f=|\nabla u|$ the following statement:

if $\g u =0$  then for all $\a>0$, it holds
\be
\lab{uweak}
\a^{n/(n-1)}  | \{ x \, | \, |u(x)|> \a \} | \le C_3 \Vert \nabla u \Vert^{n/(n-1)}_1.
\ee Here $C_3$ is a constant depending only on the controlling constants in Properties A and
B.

{\it Step 2.}
\medskip

Now we will convert the weak type inequality (\ref{uweak}) to the desired $L^1$ Sobolev inequality, using an argument based on the
 idea in \cite{FGW:1}. See also \cite{CDG:1}. Define the sets
\[
D_k = \{ x \, | \, |u(x)| > 2^k \}, \quad k \, \, \text{are integers}.
\]Then
\[
\Vert u \Vert_p =\left( \Sigma^\infty_{k=-\infty}
\int_{D_k -D_{k+1}} |u(x)|^p dg \right)^{1/p}
\] where $p=n/(n-1)$ here and later in the proof.
This shows
\be
\lab{up<sum}
\Vert u \Vert_p \le \left( \Sigma^\infty_{k=-\infty}   2^{(k+1) p} | D_k| \right)^{1/p} =
\left(  \Sigma^\infty_{k=-\infty}   2^{(k+1)p } | \{ x \, | \, |u(x)|>2^k \}| \right)^{1/p}.
\ee Now we define
\[
g_k=g_k(x)=
\begin{cases} 2^{k-1}, \qquad x \in D_k =  \{ x \, | \, |u(x)| > 2^k \},\\
|u(x)|-2^{k-1}, \qquad x \in D_{k-1}-D_k= \{ x \, | \, 2^{k-1} < |u(x)| \le 2^k \},\\
0, \qquad x \in D^c_{k-1}= \{ x \, | \, |u(x)| \le  2^{k-1} \}.
\end{cases}
\] It is clear that $g_k$ is a Lipschitz function such that $0 \le g_k \le |u|/2$.

 Observe that
\[
D_k \subset  \{ x \, | \, g_k(x) = 2^{k-1} \} \subset  \{ x \, | \, g_k(x) >2^{k-2} \}
\subset  \{ x \, | \, |g_k(x)-\g g_k| >2^{k-3} \} \cup  \{ x \, | \, \g g_k >2^{k-3} \}
\] Here $\g g_k$ is the average of $g_k$ on $\M$. Hence
\be
\lab{Dk<sum}
\al
|D_k| &\le  |\{ x \, | \, |g_k(x)-\g g_k| >2^{k-3} \}| + |  \{ x \, | \, \g g_k >2^{k-3} \}|\\
&\equiv T_{k1} + T_{k2}.
\eal
\ee Note the average of the function $g_k-\g g_k$ is $0$. Thus we can apply
(\ref{uweak}), with $u$ there being replaced by $g_k-\g g_k$, to deduce
\be
\lab{Tk1}
T_{k1}=|\{ x \, | \, |g_k(x)-\g g_k| >2^{k-3} \}| \le C C_3  2^{-p k} \Vert \nabla g_k \Vert^{p}_1.
\ee To treat $T_{k2}$, recall that $g_k \le |u|/2$ which implies
\[
\g g_k \le  \Vert u \Vert_1 /(2 |\M|).
\] Therefore
\[
T_{k2} = |  \{ x \, | \, \g g_k >2^{k-3} \}| \le  |  \{ x \, | \, \frac{\Vert u \Vert_1}{|\M|} >2^{k-2} \}|.
 \] This shows that
 \be
 \lab{Tk2}
 \al
 T_{k2} =
 \begin{cases}
 0, \quad \text{when} \quad k > 2 + \log_2 \frac{\Vert u \Vert_1}{|\M|} \equiv k_0\\
 |\M|, \quad  \quad k \le k_0.
 \end{cases}
 \eal
 \ee Substituting (\ref{Tk1}) and (\ref{Tk2}) into (\ref{Dk<sum}), we deduce
\[
\al
|D_k|
\le
 \begin{cases}
 C C_3  2^{-p k} \Vert \nabla g_k \Vert^{p}_1, \quad \text{when} \quad k > k_0\\
 C C_3  2^{-p k} \Vert \nabla g_k \Vert^{p}_1+ |\M|, \quad  \quad k \le k_0.
 \end{cases}
 \eal
 \]Substituting this to (\ref{up<sum}) and using Minkowski inequality, we obtain
 \[
 \Vert u \Vert_p \le C_4 \Sigma^\infty_{k=-\infty}  \Vert \nabla g_k \Vert_1
 + C |\M|^{1/p}  \Sigma^{[k_0] +1}_{k=-\infty} 2^k.
 \] Here $[k_0]$ is the greatest integer less than or equal to $k_0$. Note that the supports of $\nabla g_k$ are disjoint and $\nabla g_k = \nabla |u|$ in the supports.  Also by the definition of  $k_0$ in (\ref{Tk2}), we have $2^{k_0} = 4 \Vert u \Vert_1/|\M|$. Hence
 \[
  \Vert u \Vert_p \le C_4  \Vert \nabla u \Vert_1 +  C |\M|^{1/p}  \Vert u \Vert_1/|\M|,
\] which implies, since $p=n/(n-1)$,
\[
 \Vert u \Vert_{n/(n-1)} \le C_4  \Vert \nabla u \Vert_1 +  C \frac{1}{|\M|^{1/n}}  \Vert u \Vert_1.
\]Here $C$ is a numerical constant. This proves Theorem \ref{thm1.1}. \qed

Two final remarks are in order.
Let $\alpha$ be the average of $u$ in $\M$. Then By Theorem \ref{thm1.1}, we
have
\[
 \Vert u-\alpha \Vert_{n/(n-1)} \le C  \Vert \nabla u \Vert_1 +  C \frac{1}{|\M|^{1/n}}  \Vert u -\alpha \Vert_1.
\]Since the average of $u-\alpha$ is zero, inequality (\ref{u<I1du}) implies
\[
|u(x)-\alpha| \le C_0 \int_\M  \frac{|\nabla u(y)|}{d(x, y)^{n-1}} dg(y).
\]After integration, using the $\kappa$ noninflating property, we find that
\[
 \Vert u -\alpha \Vert_1 \le C diam (\M) \Vert \nabla u \Vert_1.
\]By the $\kappa$ noncollapsing property of $\M$, there holds
$diam (\M) \le C |\M|^{1/n}$. This shows the usual isoperimetric inequality
\[
 \Vert u-\alpha \Vert_{n/(n-1)} \le C  \Vert \nabla u \Vert_1.
\]

Notice the $L^2$ Poincar\'e inequality (\ref{goodl2Poin}),
\cite{Z11:1} and Section 9 of \cite{Ch:1} imply the following long
time convergence result: {\it the K\"ahler Ricci flow in Theorem
\ref{thm1.1} converges sequentially in time, under
Gromov-Hausdorff  topology, to a compact metric space with $L^2$
Poincar\'e inequality and volume doubling condition.} By Cheeger's
Theorem 11.7 \cite{Ch:1}, the limit space can be equipped with
differential structure a.e..

{\bf Acknowledgment.} Q. S. Z. would like to thank Professors L.
Capogna, X. X. Chen, D. Jerison, Bing Wang and Bun Wong for
helpful conversations. Part of the work was done when he was a
visiting professor at Nanjing University under a Siyuan Foundation
grant, the support of which is gratefully acknowledged.

Both of us wish to thank the referee for checking the paper
carefully and making helpful corrections and suggestions.

\bigskip

\noindent e-mail:  tian@math.princeton.edu and qizhang@math.ucr.edu

\enddocument